\documentclass[12pt]{article} \usepackage{latexsym}
\usepackage{amsmath,amssymb,amsthm}
\usepackage{graphicx,color,url,paralist} 
\usepackage{titling,url}

\begin{document}

\newtheorem{theorem}{Theorem}[section]
\newtheorem{lemma}[theorem]{Lemma}
\newtheorem{corollary}[theorem]{Corollary}
\newtheorem{proposition}[theorem]{Proposition}
\newtheorem{question}[theorem]{Question}
\theoremstyle{definition}
\newtheorem{definition}[theorem]{Definition}
\newtheorem{example}[theorem]{Example}
\theoremstyle{remark}
\newtheorem*{remark}{Remark}

\newcommand\RE{{\sc Random Edge}}
\newcommand\R{{\mathbb R}}
\newcommand\N{{\mathbb N}}
\newcommand{\reworst}{f}
 
\title{The \RE\ Rule\\
  on Three-Dimensional Linear Programs%
\thanks{%
  Work on this paper by Micha Sharir was supported by NSF Grants
  CCR-97-32101 and CCR-00-98246, by a grant from the U.S.-Israeli
  Binational Science Foundation, by a grant from the Israel Science
  Fund (for a Center of Excellence in Geometric Computing), and by the
  Hermann Minkowski--MINERVA Center for Geometry at Tel Aviv
  University.  Volker Kaibel and G\"unter M.~Ziegler  were
  supported by the DFG-Forschergruppe \emph{Algorithmen, Struktur,
    Zufall} (FOR 413/1-1, Zi 475/3-1).  G\"unter M.~Ziegler was
  supported by a DFG Leibniz grant.  Part of the work was done during
  the workshop ``Towards the Peak'' at La Claustra, Switzerland,
  August 2001.}\\%
(extended abstract) 
}

\thanksmarkseries{arabic}

\author{%
Volker Kaibel\thanks{%
  MA 6--2, TU Berlin, 10623~Berlin, Germany, \texttt{\{kaibel,mechtel,ziegler\}@math.tu-berlin.de}
}
\and 
Rafael Mechtel\thanksmark{2}%
\and
Micha Sharir\thanks{%
  School of Computer Science, Tel Aviv University, Tel-Aviv
  69978, Israel and Courant Institute of Mathematical Sciences, 
  New York University, New York, NY 10012, USA,
  \texttt{michas@post.tau.ac.il} 
}
\and 
G\"unter M. Ziegler\thanksmark{2}%
}

\maketitle

\begin{abstract}
The worst-case expected 
length of the path taken by the simplex algorithm
with the \RE\ pivot rule on a $3$-dimensional linear program 
with $n$ constraints is shown to be bounded by
\[
  1.3445\cdot n\ \le \ f(n)\  \le\ 1.4943\cdot n
\]
for large enough~$n$.
\end{abstract}

\section{Introduction}

The \RE\ pivot rule is undoubtedly the most natural, and simplest
(randomized) pivot rule for the simplex algorithm: ``At each iteration,
proceed from the current vertex of the polyhedron~$P$ of feasible
solutions to an improving neighbor, chosen uniformly at random in the
one-dimensional skeleton (i.e., the graph) of~$P$.''

Despite its simplicity, this algorithm until now has resisted almost
all attempts to analyze its worst-case behavior, with a few
exceptions for special cases, among them the linear assignment problem
(Tovey~\cite{Tovey}), certain linear programs on cubes, including the
Klee-Minty cubes (Kelly~\cite{Kelly}, G\"artner, Henk \&\ 
Ziegler~\cite{Z38a}), and $d$-dimensional linear programs (i.e.,
$\dim(P)=d$) with at most $d+2$ constraints (G\"artner et
al.~\cite{gaertner01:_one}). All known results leave open the
possibility that the \RE\ pivot rule yields a strongly polynomial time
algorithm -- it might be even quadratic. In particular, it is not
fooled by the deformed products (defined by Amenta and Ziegler
\cite{Z51a}), which yield the well-known exponential examples for all
the classical deterministic pivot rules.


Here, we only treat the case of $3$-dimensional linear programs,
which, of course, is solved in linear time by every (even
deterministic) finite variant of the simplex algorithm. Nevertheless,
due to the remarks made above, it seems interesting to analyze the
\RE\ pivot rule for this case -- and here too, it seems that
accurate analysis of \RE\ is quite hard.

With the usual reductions (see, e.g., Ziegler \cite[Lect.~3]{Z35}), we
may assume that our linear program is $\min\{x_3\ :\ x\in P\}$, where
$P$ is a $3$-dimensional simple polytope (its graph is $3$-regular)
with exactly $n$ facets, and hence $2n-4$ vertices and $3n-6$ edges,
and no two vertices have the same objective function value $x_3$.
Thus we have an ordering of the vertices
$v_{2n-5},v_{2n-6},\dots,v_1,v_0$ by decreasing objective function
(i.e., by height). Here $v_0=v_{\min}$ is the unique minimal (lowest)
vertex of the linear program, while $v_{2n-5}=v_{\max}$ is the unique
maximal (highest) vertex.

The expected 
length of the path (i.e., the number of pivot steps) taken by the
simplex algorithm
on the linear program, starting at vertex $v$ of $P$ and using the \RE\ rule, 
is then given by $E(v_0)=0$ and
$$
E(v)\ =\ 1 + \frac{1}{d_v}\sum_{j=1}^{d_v}E(w_j)\qquad(v\not= v_0)\ ,
$$
where $w_1,\dots,w_{d_v}$ are the lower neighbors of~$v$. It is
easy to see that (in addition to the unique maximal vertex and 
the unique minimal vertex) there are $n-3$ vertices~$v$ with $d_v=1$
(\emph{$1$-vertices}) and $n-3$ vertices~$v$ with $d_v=2$
(\emph{$2$-vertices}); this is the $3$-dimensional case of the
Dehn-Sommerville equations \cite[Thm.~8.21]{Z35}.

Define $\reworst(n)$ to be the maximum expected number of pivot steps
taken by the \RE\ algorithm on any $3$-dimensional linear program with
$n$ constraints. 
While it is quite straightforward to construct a sequence of examples
with $E(v)\ge\frac43\cdot n - \text{const}$, our results in
Section~\ref{sec:low} (Theorem~\ref{thm:low}) show
$$
f(n)\ \ge \ \frac{1721}{1280}\cdot n - \frac{4722}{1280}
$$
for infinitely many~$n$ in arithmetic progression
($\frac{1721}{1280}=\frac43+\frac{43}{3840}$).
In Section~\ref{sec:high} we prove
$$
f(n)\ \le \ \frac{130}{87}\cdot n -\frac{115}{29}
$$
(Theorem~\ref{up:main}).  Both results taken together, this yields
that
$$
  1.3445\cdot n\ \le \ f(n)\  \le\ 1.4943\cdot n
$$
holds for all large enough~$n$. In particular, asymptotically $f(n)$
lies between $(\frac43+\varepsilon)\cdot n$ and
$(\frac32-\varepsilon)\cdot n$ for some $\varepsilon>0$.
Determining
``the right coefficient''
seems, however, to be very hard.

\section{Lower Bounds}
\label{sec:low}

The expected 
number of pivot steps required by
the simplex algorithm using \RE\ only depends
on the graph of the polytope, directed via the objective function.
Therefore, we will describe our examples yielding lower bounds on
$f(n)$ by the corresponding directed graphs. The following result
provides a nice certificate for a directed graph to come from a
$3$-dimensional linear program. 

\begin{theorem}[Mihalisin and Klee \cite{MihalisinKlee}]
  A directed graph~$D$ (without loops and parallel arcs) is induced by a
  $3$-dimensional linear program if and only if 
  \begin{compactenum}
  \item[$\circ$] it is planar and $3$-connected (as an undirected graph),
  \item[$\circ$] it is acyclic with a unique source and a unique sink,
  \item[$\circ$] it has a unique local sink in every face cycle (these are the
    non-separating induced cycles), and
  \item[$\circ$] it admits three directed paths from its source to
        its sink that have disjoint sets of interior nodes.
  \end{compactenum}
\end{theorem}

\subsection{Duals of Cyclic Polytopes}
\label{subsec:dualcyc}

\begin{figure}[ht]
  \centering
  \input{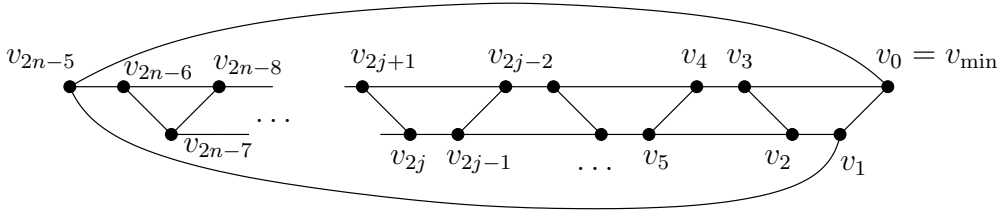}  
  \caption{The example on the dual cyclic polytope.}
  \label{fig:dualcyc}
\end{figure}

\paragraph{Example~1.}
Our first sequence of examples are wedges, i.e., they are
combinatorially equivalent to duals of cyclic polytopes.
Figure~\ref{fig:dualcyc} depicts the orientations of the edges. Here,
as well as in the sequel, our convention is that the ordering of the
vertices from left to right in the figure defines the (decreasing)
ordering of the vertices according to the objective function. It is
easy to see that the conditions of the Mihalisin--Klee theorem are
satisfied.

For the expected 
number of pivot steps
$E(v_i)$, we then have the starting
values $E(v_0)=0$ and $E(v_1)=1$, and the recurrences
\[
E(v_{2j})\ =\ E(v_{2j-1})+1,\qquad
E(v_{2j+1})=\tfrac12\big(E(v_{2j})+E(v_{2j-2})\big)+1
\]
for $0<j\le n-4$. Thus, using induction, we obtain 
\[
E(v_{2j}) + 2 E(v_{2j+1})\ =\ 4j+2
\]
for $0\le j\le n-4$.
In particular, for $j=n-4$ this yields
\[
\max\Big\{E(v_{2n-8}),E(v_{2n-7})\Big\}\ \ge\ \frac{4n}{3}-\frac{14}{3}.
\]

\subsection{Improved Lower Bounds}
\label{subsec:improved}

The next examples are based on the construction of a ``backbone
polytope'': This will be a simple $3$-polytope $P_k$ with $k+2$ facets
and $2k$ vertices, of which $k$ vertices $v_{k-1},\dots,v_0$ form a
decreasing chain, such that $v_0$ is the minimal vertex, and $v_{i-1}$
is the only lower neighbor of $v_i$, for $i>0$.

\paragraph{Constructing the backbone.} 
We start with the simplex, obtained for $k=2$ with vertices $w_0,w_1,v_1,v_0$,
as shown on the left. We then inductively cut off the vertex $v_{k-1}$ by a plane,
replacing it by a small triangle, as shown on the right.
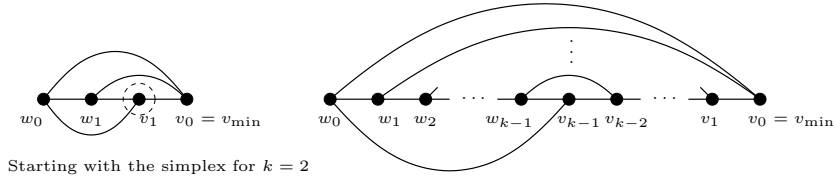
\begin{center}
\input{backbone2.pstex_t}
\end{center}

\paragraph{Example 2.}
Our second sequence of examples is obtained from the backbone
polytopes $P_k$ by performing three specific vertex cuts at each
vertex $v_i$, for $i=0, \ldots, k-1$.  Before cutting (in $P_k$), each
vertex $v_i$ ($i>0$) has indegree~$2$, while $v_0$ has indegree $3$.
The two vertex cuts are supposed to create at each vertex $v_i$ the
following configuration (where again all edges are directed ``to the
right''):

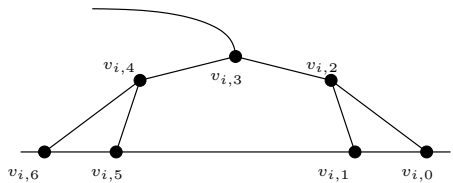
\begin{center}
\input{template2.pstex_t}
\end{center}
This creates a simple polytope $P_k'$ with $n=k+2+3k=4k+2$ facets.

Our starting vertex for \RE\ on this example will be $v_{k-1,6}$; the
expected number of steps taken by \RE\ is the sum of the
probabilities $p_e$ that the edge $e$ is traversed.  We think of these
probabilities as a ``flow'' from $v_{k-1,6}$ to $v_{0,0}$. Our figure
indicates the flow values on the edges, for a flow of total value $8$;
equivalently, these are the transversal probabilities in units
of~$\frac18$.

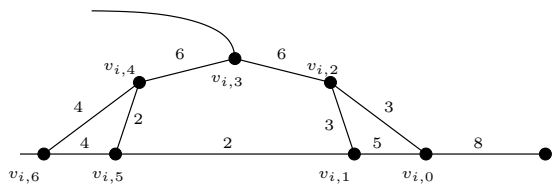
\begin{center}
\input{template2_ana.pstex_t}
\end{center}

We get the same values for each of the triple-vertex-cut-off
configurations, except for the last one, which has no edge leaving the
global sink $v_{0,0}$.  Thus the expected number of \RE\ steps,
starting from $v_{k-1,6}$, is
\[
E(v_{k-1,6})\ \ =\ \ k\cdot
 \frac{43}{8} - 1 = \frac{43}{32}n-\frac{59}{16}\ ,
\]
with $k=\frac14(n-2)$. Asymptotically, this yields a better lower
bound, due to $\frac{43}{32}>\frac43$.

The graph of $P_k'$ looks like this:
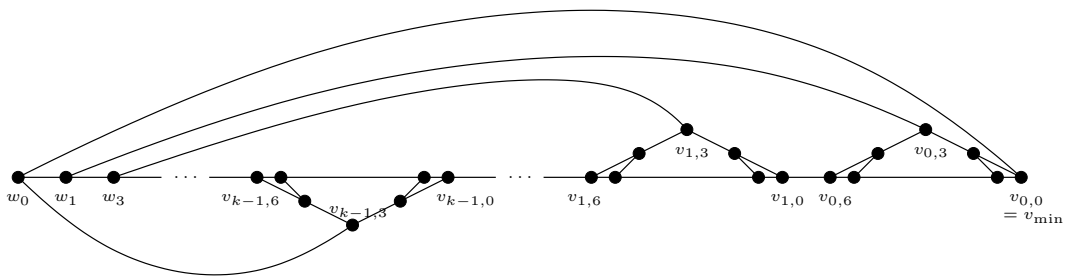
\begin{center}
\input{complete.43.32.pstex_t}
\end{center}

\paragraph{Example~3.}%
The last examples produced $k = \frac 14(n-2)$ vertices which are not
used.  We will further improve the lower bound by using more facets in
each local configuration and thus reducing the number of unused
vertices (though our final example will still have linearly many
unused vertices). 

We use the same backbone polytope as before, but we replace each
vertex $v_i$, for $i=0, \ldots, k-1$, with the following graph.  (We
do not give an explicit polytopal construction for this example, but
it can be constructed by cutting off vertices and edges of the
backbone polytope. Alternatively, such a construction is provided by
the Mihalisin-Klee Theorem.)  To analyze the random path length
through this graph we send from $v_{i,18}$ $128$ units of flow through
the network.

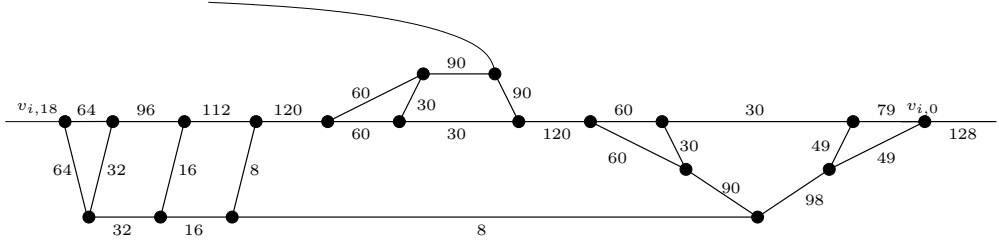
\begin{center}
\input{theworst2.pstex_t}
\end{center}
The total flow through all edges is $1721$. For each configuration $9$
facets are required. Together with the facets from the backbone construction,
this yields $n = 10k+2$. Hence we get
\[
E(v_{k-1,18})\ =\ k\cdot\frac{1721}{128}-1\ =\ 
\frac{1721}{1280}n - \frac{4722}{1280}\ ,
\]
where
$\frac{1721}{1280} = \frac{43}{32}+\frac{1}{1280} = 
\frac 43 + \frac{1}{96} + \frac{1}{1280} > 1.3445$.

In contrast to the preceeding examples, this example does not contain a
directed Hamiltonian path.

Summarizing, we have proved the following bound.

\begin{theorem}
\label{thm:low}
  For $n=10k+2\ge 12$,
  $$
  f(n)\ge\frac{1721}{1280}n - \frac{4722}{1280}
  $$
  holds.
\end{theorem}

\paragraph{Starting from the source.}
By splitting the maximal vertex, one can also construct examples where
the expected number of steps \emph{starting at the maximal vertex} is
at least $(\frac{1721}{1280}-\varepsilon)n$. This observation is due
to G\"unter Rote.

\section{Upper Bounds}
\label{sec:high}

Consider any linear program on a simple $3$-polytope with the
notations as described in the introduction.  For a vertex $v$, let
$N_1(v)$ (resp., $N_2(v)$) denote the number of 1-vertices (resp.,
2-vertices) that are not higher than $v$ (including $v$ itself).
Put $N(v)=N_1(v)+N_2(v)$. This is the number of vertices lower than
$v$. The core of our upper bound on $f(n)$ is the following result.

\begin{theorem} \label{32}
  For each vertex $v$, other than the maximal vertex~$v_{2n-5}$,
  we have
$$
E(v) \le \frac{46}{87} N_1(v)+ \frac{42}{87} N(v)\ .
$$
\end{theorem}
Theorem~\ref{32} implies
$$
E(v)\le \frac{130(n-3)}{87}+\frac{15}{29}=\frac{130}{87}\cdot n -\frac{115}{29} 
$$
for all~$v$. Here, the ``$\frac{15}{29}$'' comes from the fact that
Theorem~\ref{32} is proved only for $v\not=v_{2n-5}$; therefore, we
bound $E(v_{2n-5})$ by $1 + \frac{1}{3}\sum_{i=1}^3 E(w_i)$, where
$w_1,w_2,w_3$ are the neighbors of~$v_{2n-5}$, and we exploit
$N(w_1)+N(w_2)+N(w_3)\le 6n-21$.

\begin{theorem} \label{up:main}
For every 3-dimensional linear program with $n$ constraints, the expected
number of pivot steps taken by \RE\ is not more than
$$
\frac{130}{87}\cdot n -\frac{115}{29}\ .
$$
\end{theorem} 

In the remaining part of this section, we briefly sketch the proof of
Theorem~\ref{32}.  It proceeds by deriving the generic inequality
\begin{equation}
  \label{eq:generic}
  E(v) \le \alpha N_1(v)+ \beta N(v)\ ,
\end{equation}
where, for most of the proof, $\alpha$ and $\beta$ are treated as
indeterminates. Each step of the proof yields a linear inequality on
$\alpha$ and $\beta$ that needs to be satisfied in order to
imply~(\ref{eq:generic}). The proof is then completed once it is shown
that $(\alpha,\beta)=(\frac{46}{87},\frac{42}{87})$ satisfies all
inequalities; in fact, it is optimal with respect to the objective
function $\alpha+2\beta$ (see below for more details).

Inequality~(\ref{eq:generic}) is proved by induction on $N(v)$. The
base case $N(v)=0$ is obvious, since $v$ is the optimum in this case,
and $E(v)=0$.  Suppose the theorem holds for all vertices lower than
some vertex $v$.

We express $E(v)$ in terms of the expected costs $E(w)$ of certain
vertices $w$ that are reachable from $v$ via a few downward edges. The
general form of such a recursive expression will be
$$
E(v) = c + \sum_{i=1}^k \lambda_i E({w_i})\ ,
$$
where $\lambda_i > 0$ for each $i=1,\dots,k$,
and $\sum_i\lambda_i = 1$.

Since we assume by induction that 
$E({w_i}) \le \alpha N_1(w_i) + \beta N(w_i)$, 
for each $i$, it suffices to show that
$$
\sum_{i=1}^k \alpha \lambda_i \left(N_1(v)-N_1(w_i)\right) +
\sum_{i=1}^k \beta \lambda_i \left(N(v)-N(w_i)\right) \ge c\ .
$$
Write
$$
\Delta_1(w_i) = N_1(v)-N_1(w_i)\ ,
$$
$$
\Delta(w_i) = N(v)-N(w_i)\ ,
$$
for $i=1,\ldots,k$.
(Of course, these terms are defined with respect to the 
currently considered
vertex $v$.)
Note that $\Delta(w_i)$ is the 
\emph{distance} between~$v$ and~$w_i$, 
that is,
one plus the number of vertices between $v$ and $w_i$.

We thus need to show that
\begin{equation} \label{ind}
\sum_{i=1}^k \alpha \lambda_i \Delta_1(w_i) +
\sum_{i=1}^k \beta \lambda_i \Delta(w_i) \ge c\ .
\end{equation} 
This requires a quite extensive case analysis, of which we present
only the beginning in this extended abstract in order to indicate the
kinds of arguments used. The complete case analysis is given in the appendix.

\paragraph{Case 1:} $v$ is a 1-vertex.\\ 
Let $w_1$ denote the target of the unique downward edge emanating from
$v$ as in the following figure, where (here and in all subsequent
figures) each edge is labelled by the
probability of reaching it from $v$.
\begin{center}
\input redge1.pstex_t
\end{center}
In this case, $E(v) = 1+E(w_1)$ holds. In the setup presented above,
we have
$c=1$, $\Delta_1(w_1) \ge 1$, and $\Delta(w_1) \ge 1$,
thus (\ref{ind}) is implied by
\begin{equation} \label{eq1}
\alpha+\beta\ge 1\ .
\end{equation} 

\paragraph{Case 2:} $v$ is a 2-vertex.\\
Let $w_1$ and $w_2$ denote the targets of the two downward edges
emanating from $v$, where $w_2$ is lower than $w_1$.
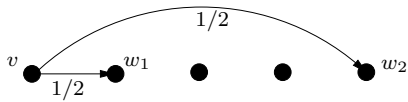
\begin{center}
\input redge2.pstex_t
\end{center}
We have
$$
E(v) = 1+\frac{1}{2}E({w_1}) + \frac{1}{2}E({w_2})\ ,
$$
hence we need to require that
$$
\frac{\alpha}{2}\Delta_1(w_1) + 
\frac{\alpha}{2}\Delta_1(w_2) + 
\frac{\beta}{2}\Delta(w_1) + 
\frac{\beta}{2}\Delta(w_2) \ge 1\ .
$$
Note that $\Delta(w_1)\ge 1$.

\paragraph{Case 2.a:} $\Delta(w_2)\ge 4$.\\
Ignoring the effect of the $\Delta_1(w_j)$'s, it suffices to require
that
$$
\frac{\beta}{2}\Delta(w_1) + \frac{\beta}{2}\Delta(w_2) \ge 1\ ,
$$
which will follow if
\begin{equation} \label{eq2}
\beta\ge \frac{2}{5}\ .
\end{equation}

\paragraph{Case 2.b.i:} $\Delta(w_2) = 3$ and one of the two vertices above
$w_2$ and below $v$ is a 1-vertex.\\
In this case $\Delta_1(w_2)\ge 1$ and $\Delta(w_1)+\Delta(w_2)\ge 4$,
so (\ref{ind}) is implied by
\begin{equation} \label{eq3}
\frac{1}{2}\alpha+2\beta\ge 1\ .
\end{equation}

We skip the remaining cases~2.b.ii and~2.c ($\Delta(w_2) = 2$) in this
extended abstract; the latter one splits into quite a large number of
subcases, which become slightly more involved. One ends up with
roughly 24 linear inequalities in addition to~(\ref{eq1}),
(\ref{eq2}), (\ref{eq3}).

Assuming we have $\alpha$ and $\beta$ that satisfy all these
inequalities, each of the induction steps is justified, and the
inequality~(\ref{eq:generic})
follows. Since we always have
$N_1(v),N_2(v) \le n-3$, we obtain
$$
E(v) \le \alpha N_1(v) + \beta(N_1(v)+N_2(v)) =
(\alpha+\beta) N_1(v) + \beta N_2(v) \le (\alpha+2\beta)(n-3)\ .
$$
Hence we choose $(\alpha,\beta)$ to minimize $\alpha+2\beta$,
subject to all the derived inequalities. This is indeed the choice
appearing in the statement of Theorem~\ref{32}.

\section{Discussion}
The improved lower bounds of Section~\ref{subsec:improved} arose
from complete enumerations for small~$n$. In particular, the lower
bounds provided by examples~2 and~3 are tight for $n=10,12$, respectively. Thus, we have
$f(10)=\frac{39}{4}=9.75$
and $f(12)=\frac{1593}{128}\approx 12.45$.

We are convinced that the bound in Theorem~\ref{up:main} is not tight.
In fact, precisely two of the inequalities of the proof of
Theorem~\ref{32} are tight for
$(\alpha,\beta)=(\frac{46}{87},\frac{42}{87})$.  The two corresponding
subcases thus constitute the bottleneck for the current upper bound.
In order to improve the bound, one should expand these two subcases
further, aiming at replacing those two inequalities by weaker ones (at
the cost of a longer proof).  As a matter of fact, in an earlier
(unpublished) version of this manuscript we had obtained an upper
bound of $1.5\cdot n$, using a somewhat more compact enumeration
scheme. The current scheme is a refinement, based on further expansion
of the preceding one.

At this point, we have no real sense of what the exact bound should
be. The refinement of the approach in this section, as just
outlined, is not likely to yield substantial improvements in the upper
bound, so a radically different approach is probably called for.
Such an improvement might be based on the observation that
certain local structures involve 1-vertices with one of its
upward neighbors lying above $v$. In fact, if the portion below
$v$ contains $k$ such vertices, there must exist at least $k$
vertices higher than $v$, so the upper bound 
$\alpha N_1(v) + \beta N(v)$ is much smaller than
$(\alpha+2\beta)(n-3)$. As a matter of fact, the lower bounds
derived in Section~\ref{sec:low} do take this constraint into
consideration.

Another observation is that the proof of Theorem~\ref{32} (as
detailed in the appendix) uses (twice)
the 3-connectivity of the edge graph of $P$, but it does not use its
planarity at all, although it does occassionally run into nonplanar
configurations.  It is conceivable that further refinement stages
might reach nonplanar configurations, whose exclusion would allow us
to further improve the bound.

What if we also drop the 3-connectivity assumption? Then we need to
consider additional cases, which cause our upper bound to increase.
The best upper bound we have at the moment for this relaxed situation
is $13n/8=1.625\cdot n$, but we are convinced that it too can be
further improved.

\paragraph{Acknowledgements.}
We are grateful to Emo Welzl and G\"unter Rote for inspiring
discussions and helpful comments.

\bibliographystyle{siam}
\bibliography{random3d}

\clearpage
\section*{Appendix}
This appendix contains the complete case analysis of the proof of
Theorem~3.1 in the extended abstract (including the cases presented
there).  

\setcounter{equation}{2}
\paragraph{Case 1:} $v$ is a 1-vertex.\\ 
Let $w_1$ denote the target of the unique downward edge emanating from
$v$ as in the following figure, where each edge is labelled by the
probability of reaching it from $v$.
\begin{center}
\input redge1.pstex_t
\end{center}
In this case, $E(v) = 1+E(w_1)$ holds. In the setup presented above,
we have
$c=1$, $\Delta_1(w_1) \ge 1$, and $\Delta(w_1) \ge 1$,
thus (2) is implied by
\begin{equation} \label{app:eq1}
\alpha+\beta\ge 1\ .
\end{equation} 

\paragraph{Case 2:} $v$ is a 2-vertex.\\
Let $w_1$ and $w_2$ denote the targets of the two downward edges
emanating from $v$, where $w_2$ is lower than $w_1$.
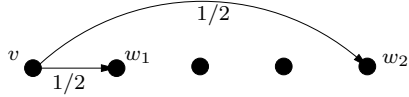
\begin{center}
\input redge2.pstex_t
\end{center}
We have
$$
E(v) = 1+\frac{1}{2}E({w_1}) + \frac{1}{2}E({w_2})\ ,
$$
hence we need to require that
$$
\frac{\alpha}{2}\Delta_1(w_1) + 
\frac{\alpha}{2}\Delta_1(w_2) + 
\frac{\beta}{2}\Delta(w_1) + 
\frac{\beta}{2}\Delta(w_2) \ge 1\ .
$$
Note that $\Delta(w_1)\ge 1$.

\paragraph{Case 2.a:} $\Delta(w_2)\ge 4$.\\
Ignoring
the effect of the $\Delta_1(w_j)$'s, it suffices to
require that
$$
\frac{\beta}{2}\Delta(w_1) + \frac{\beta}{2}\Delta(w_2) \ge 1\ ,
$$
which will follow if
\begin{equation} \label{app:eq2}
\beta\ge \frac{2}{5}\ .
\end{equation}

\paragraph{Case 2.b.i:} $\Delta(w_2) = 3$ and one of the two vertices above
$w_2$ and below $v$ is a 1-vertex.\\
In this case
$\Delta_1(w_2)\ge 1$
and $\Delta(w_1)+\Delta(w_2)\ge 4$, so (2) is implied by
\begin{equation} \label{app:eq3}
\frac{1}{2}\alpha+2\beta\ge 1\ .
\end{equation}

\paragraph{Case 2.b.ii:} $\Delta(w_2) = 3$ and the two vertices between
$v$ and $w_2$ are 2-vertices.
Denote the second intermediate vertex as $v'$.
We may assume that $v'$ is reachable from $v$,
otherwise we can ignore it and reduce the situation to Case 2.c
treated below
(be choosing another ordering of the vertices producing the same
oriented graph)%
. Three subcases can arise.

First, assume that none of the three edges that emanate from $w_1$ and
$v'$ further down reaches $w_2$. Denote by $x,y$ the two downward
neighbors of $v'$ and by $z$ the downward neighbor of $w_1$ other than
$v'$. The vertices $x,y,z$ need not be distinct but none of them
coinicdes with $w_2$.
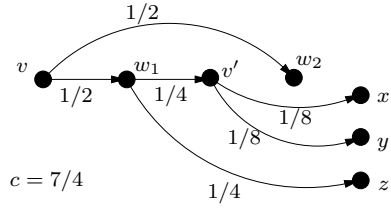
\begin{center}
\input redge2bii1.pstex_t
\end{center}
We have here $c=7/4$.

To make the analysis simpler to follow visually, we present it in a
table. Each row denotes one of the target vertices $w_2,x,y,z$,
`multiplied' by the probability of reaching it from $v$. The left 
(resp., right) column denotes a lower bound on the corresponding quantities 
$\Delta_1(\cdot)$ (resp., $\Delta(\cdot)$). 
To obtain an inequality that implies (2), one has to multiply
each entry in the left (resp., right) column by the row probability
times $\alpha$ (resp., times $\beta$), and require that the sum of all
these terms be $\ge c$.
\begin{center}
\begin{tabular}{|l|c|c|}
& $\alpha \Delta_1$ & $\beta \Delta$ \\
\hline
$1/2 w_2$ & 0 & 3 \\
$1/8 x$ & 0 & 4 \\
$1/8 y$ & 0 & 5 \\
$1/4 z$ & 0 & 4 \\
\hline
\end{tabular}
\end{center}
Note the following: (a) We do not assume that the rows represent 
distinct vertices (in fact, $x=z$ is implicit in the table); 
this does not cause any problem
in applying the rule for deriving an inequality from the table.
(b) We 
have
to squeeze the vertices so as to make the resulting
inequality as sharp (and difficult to satisfy) as possible; 
thus we made one of $x,y$ the farthest
vertex, because making $z$ the farthest vertex would have made the
inequality easier to satisfy.

We thus obtain
$$
\left(\frac{3}{2}+\frac{4}{8}+\frac{5}{8}+\frac{4}{4}\right)\beta\ge
\frac{7}{4} ,
$$
or
\begin{equation} \label{app:eq4}
\beta\ge \frac{14}{29} .
\end{equation}

Next, assume that $w_2$ is 
connected to
$v'$. In this case $w_2$ is a
1-vertex, and we extend the configuration to include its unique downward
neighbor $w_3$.
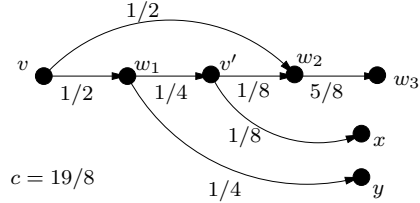
\begin{center}
\input redge2bii2.pstex_t
\end{center}
Let $x$ denote the other downward neighbor of $v'$ and 
let $y$ denote the other downward neighbor of $w_1$. 
In the following table, the `worst' case is to make $w_3$ and $y$ 
coincide, and make $x$ the farthest vertex.
\begin{center}
\begin{tabular}{|l|c|c|}
& $\alpha \Delta_1$ & $\beta \Delta$ \\
\hline
$5/8 w_3$ & 1 & 4 \\
$1/8 x$ & 1 & 5 \\
$1/4 y$ & 1 & 4 \\
\hline
\end{tabular}
\end{center}
We then obtain
$$
\alpha + \left(\frac{20}{8}+\frac{5}{8}+\frac{4}{4}\right)\beta\ge
\frac{19}{8} ,
$$
or
\begin{equation} \label{app:eq5}
\alpha + \frac{33}{8}\beta\ge \frac{19}{8} .
\end{equation}

Finally, assume that $w_2$ is 
connected to
$w_1$. Here too $w_2$ is a
1-vertex, and we extend the configuration to include its unique downward
neighbor $w_3$. 
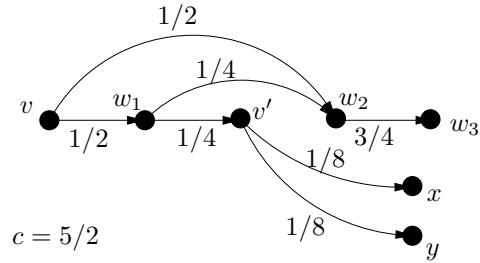
\begin{center}
\input redge2bii3.pstex_t
\end{center}
Denoting by $x,y$ the two downward neighbors of $v'$, our 
table and resulting inequality become
\begin{equation} \label{app:eq6}
\begin{tabular}{|l|c|c|}
& $\alpha \Delta_1$ & $\beta \Delta$ \\
\hline
$3/4 w_3$ & 1 & 4 \\
$1/8 x$ & 1 & 4 \\
$1/8 y$ & 1 & 5 \\
\hline
\end{tabular}
\hspace{2cm}
\alpha + \frac{33}{8}\beta\ge \frac{5}{2} ,
\end{equation}
which, by the way, is stronger than (\ref{app:eq5}).

\paragraph{Case 2.c:} $\Delta(w_2) = 2$.
Hence, the only remaining case is that $w_1$ and $w_2$ are the
two vertices immediately following $v$.

\paragraph{Case 2.c.i:} $w_1$ is a 1-vertex (whose other upward neighbor
lies above $v$). Its unique downward edge ends at some vertex 
which is either $w_2$ or lies below $w_2$. 

Assume first that this vertex coincides with $w_2$, which makes $w_2$ a
1-vertex, whose unique downward neighbor is denoted as $v'$.
The local structure, table, and inequality are
\begin{equation} \label{app:eq7}
\input redge2ci1.pstex_t
%
\hspace{1.5cm}
\begin{tabular}{|l|c|c|}
& $\alpha \Delta_1$ & $\beta \Delta$ \\
\hline
$v'$ & 2 & 3 \\
\hline
\end{tabular}
\hspace{1.5cm}
2\alpha + 3\beta\ge \frac{5}{2} .
\end{equation}

Suppose next that the downward neighbor $w_3$ of $w_1$ lies below $w_2$.
We get
\begin{equation} \label{app:eq8}
\input redge2ci2.pstex_t
%
\hspace{1.5cm}
\begin{tabular}{|l|c|c|}
& $\alpha \Delta_1$ & $\beta \Delta$ \\
\hline
$1/2 w_2$ & 1 & 2 \\
$1/2 w_3$ & 1 & 3 \\
\hline
\end{tabular}
\hspace{1.5cm}
\alpha + \frac{5}{2}\beta\ge \frac{3}{2} .
\end{equation}

\paragraph{Case 2.c.ii:} $w_1$ is a 2-vertex, both of whose downward 
neighbors lie strictly below $w_2$. Denote these neighbors as $w_3,w_4$, 
with $w_3$ lying above $w_4$.
\begin{center}
\input redge2cii.pstex_t
\end{center}

We may assume $\Delta(w_3)=3$ (i.e., there is no vertex between~$w_2$ 
and~$w_3$), since
$\Delta(w_3)\ge 4$ requires $\beta\ge\frac{6}{13}$ as the  sharpest
inequalitity, which is already implied by~(\ref{app:eq4}).
%

\paragraph{Case 2.c.ii.1:} $w_2$ is a 1-vertex. 
Then 
the table and inequality become
\begin{equation} \label{app:eq9}
\begin{tabular}{|l|c|c|}
& $\alpha \Delta_1$ & $\beta \Delta$ \\
\hline
$1/2 w_2$ & 0 & 2 \\
$1/4 w_3$ & 1 & 3 \\
$1/4 w_4$ & 1 & 4 \\
\hline
\end{tabular}
\hspace{2cm}
\frac{1}{2}\alpha + \frac{11}{4}\beta\ge \frac{3}{2} .
\end{equation}

\paragraph{Case 2.c.ii.2:} $w_2$ is a 2-vertex but $w_3$ is a 1-vertex. 
Then $w_3$ (which satisfies $\Delta(w_3)=3$) is
connected either to $w_2$ or to a vertex above $v$. 
In the former
case, let $x$ denote the other downward neighbor of $w_2$, and  
let $y$ denote the unique downward neighbor of $w_3$. The local
structure looks like this (with $x,y,w_4$ not necessarily distinct):
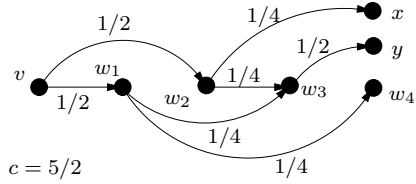
\begin{center}
\input redge2cii1.pstex_t
\end{center}
The table depends on whether $x$ precedes or succeeds
$w_3$. In the former case the (worst) table and inequality are
\begin{equation} \label{app:eq10}
\begin{tabular}{|l|c|c|}
& $\alpha \Delta_1$ & $\beta \Delta$ \\
\hline
$1/4 x$ & 0 & 3 \\
$1/2 y$ & 1 & 5 \\
$1/4 w_4$ & 1 & 5 \\
\hline
\end{tabular}
\hspace{2cm}
\frac{3}{4}\alpha + \frac{9}{2}\beta\ge \frac{5}{2} .
\end{equation}

In the latter case the (worst) table and inequality are
\begin{equation} \label{app:eq11}
\begin{tabular}{|l|c|c|}
& $\alpha \Delta_1$ & $\beta \Delta$ \\
\hline
$1/4 x$ & 1 & 4 \\
$1/2 y$ & 1 & 4 \\
$1/4 w_4$ & 1 & 5 \\
\hline
\end{tabular}
\hspace{2cm}
\alpha + \frac{17}{4}\beta\ge \frac{5}{2} .
\end{equation}

The next case is where the other upward neighbor of $w_3$ lies above
$v$. Let $x,y$ denote the two downward neighbors of $w_2$, and let $z$ 
denote the unique downward neighbor of $w_3$. The local structure is:
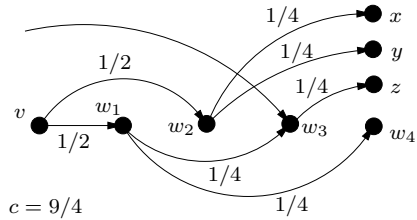
\begin{center}
\input redge2cii2.pstex_t
\end{center}
The table depends on how many of $x,y$ precede $w_3$.
If both precede $w_3$, the table and inequality become
\begin{equation} \label{app:eq12}
\begin{tabular}{|l|c|c|}
& $\alpha \Delta_1$ & $\beta \Delta$ \\
\hline
$1/4 x$ & 0 & 3 \\
$1/4 y$ & 0 & 4 \\
$1/4 z$ & 1 & 6 \\
$1/4 w_4$ & 1 & 6 \\
\hline
\end{tabular}
\hspace{2cm}
\frac{1}{2}\alpha + \frac{19}{4}\beta\ge \frac{9}{4} .
\end{equation}

If only one of $x,y$ precedes $w_3$, say $x$, the table and inequality become
\begin{equation} \label{app:eq13}
\begin{tabular}{|l|c|c|}
& $\alpha \Delta_1$ & $\beta \Delta$ \\
\hline
$1/4 x$ & 0 & 3 \\
$1/4 y$ & 1 & 5 \\
$1/4 z$ & 1 & 5 \\
$1/4 w_4$ & 1 & 6 \\
\hline
\end{tabular}
\hspace{2cm}
\frac{3}{4}\alpha + \frac{19}{4}\beta\ge \frac{9}{4} ,
\end{equation}
which is weaker than (\ref{app:eq12}).

Finally, if none of $x,y$ precedes $w_3$, the table and inequality become
\begin{equation} \label{app:eq14}
\begin{tabular}{|l|c|c|}
& $\alpha \Delta_1$ & $\beta \Delta$ \\
\hline
$1/4 x$ & 1 & 4 \\
$1/4 y$ & 1 & 4 \\
$1/4 z$ & 1 & 5 \\
$1/4 w_4$ & 1 & 5 \\
\hline
\end{tabular}
\hspace{2cm}
\alpha + \frac{9}{2}\beta\ge \frac{9}{4} .
\end{equation}

\paragraph{Case 2.c.ii.3:} Both $w_2$ and $w_3$ are 2-vertices. 
We have to consider the following type of configuration (where 
$x,y,z,t,w_4$ need not all be distinct, but $x\not= y$ and $z\not= t$, 
and we may assume $x\not=t$, $y\not= z$; also, because $\Delta(w_3)=3$,
both~$x$ and~$y$ are lower than~$w_3$):
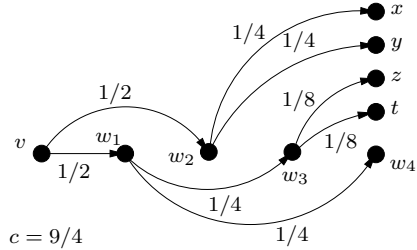
\begin{center}
\input redge2cii3.pstex_t
\end{center}
Intuitively, a worst table is obtained by `squeezing' $x,y,z,t$, and 
$w_4$ as much to the left as possible, placing two of them at distance 4
from $v$, two at distance 5, and one at distance 6. However, squeezing
them this way will make some pairs of them coincide and form 1-vertices,
which will affect the resulting tables and inequalities.

Suppose first that among the three `heavier' targets $x,y,w_4$, at most
one lies at distance 4 from $v$. The worst table and the associated
inequality are (recall that $x\ne y$):
\begin{equation} \label{app:eq14a}
\begin{tabular}{|l|c|c|}
& $\alpha \Delta_1$ & $\beta \Delta$ \\
\hline
$1/4 x$ & 0 & 4 \\
$1/4 y$ & 0 & 5 \\
$1/8 z$ & 0 & 4 \\
$1/8 t$ & 0 & 6 \\
$1/4 w_4$ & 0 & 5 \\
\hline
\end{tabular}
\hspace{2cm}
\frac{19}{4}\beta\ge \frac{9}{4} .
\end{equation}

%

Suppose then that among $\{w_4,x,y\}$, two are at distance 4 from $v$,
say $w_4$ and $y$. Then $w_4=y$ is a 1-vertex, and we denote by $w$ its
unique downward neighbor. The local structure is:
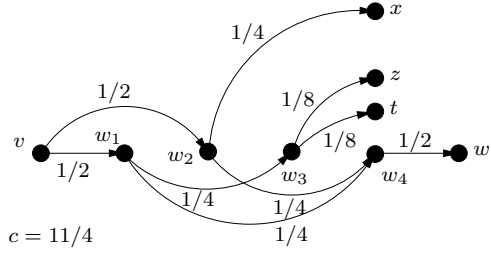
\begin{center}
  \input redge2cii4.pstex_t
\end{center}
Two equally worst tables, and the resulting common inequality are

\begin{equation} \label{app:eq:new2}
\begin{tabular}{|l|c|c|}
& $\alpha \Delta_1$ & $\beta \Delta$ \\
\hline
$1/4 x$   & 1 & 5 \\
$1/8 z$   & 1 & 6 \\
$1/8 t$   & 1 & 7 \\
$1/2 w$ & 1 & 5 \\
\hline
\end{tabular}
\hspace{1.5cm}
\begin{tabular}{|l|c|c|}
& $\alpha \Delta_1$ & $\beta \Delta$ \\
\hline
$1/4 x$   & 1 & 6 \\
$1/8 z$   & 1 & 6 \\
$1/8 t$   & 1 & 5 \\
$1/2 w$ & 1 & 5 \\
\hline
\end{tabular}
\hspace{1.5cm}
  \alpha+\frac{43}{8}\beta\ge \frac{11}{4}\ .
\end{equation}

%
%

\paragraph{Case 2.c.iii:} $w_1$ is a 2-vertex that reaches $w_2$. Then
$w_2$ is a 1-vertex, and we denote by $x$ its unique downward neighbor.
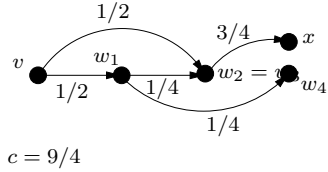
\begin{center}
\input redge2ciii.pstex_t
\end{center}
A crucial observation is that $x$ cannot be equal to $w_4$. Indeed, if 
they were equal, then $w_4$ would be a 1-vertex. 
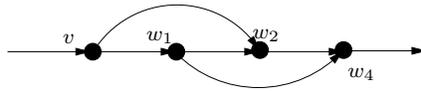
\begin{center}
\input redge2ciiix.pstex_t
\end{center}
In this case, cutting the edge graph $G$ of $P$
at the downward edge emanating from $x$ and at the edge entering $v$
would have disconnected $G$, contradicting the fact that $G$
is 3-connected%
.

We first dispose of the case where $x$ lies lower than $w_4$. 
The table and inequality are
\begin{equation} \label{app:eq15}
\begin{tabular}{|l|c|c|}
& $\alpha \Delta_1$ & $\beta \Delta$ \\
\hline
$3/4 x$ & 1 & 4 \\
$1/4 w_4$ & 1 & 3 \\
\hline
\end{tabular}
\hspace{2cm}
\alpha + \frac{15}{4}\beta\ge \frac{9}{4} .
\end{equation}
In what follows we thus assume that $x$ lies above $w_4$.

\paragraph{Case 2.c.iii.1:} $x$ is a 1-vertex that precedes $w_4$.
Suppose first that $w_4$ is the unique downward neighbor of $x$. Then
$w_4$ is a 1-vertex, and we denote its unique downward neighbor by $z$.
The local structure, table and inequality are: 
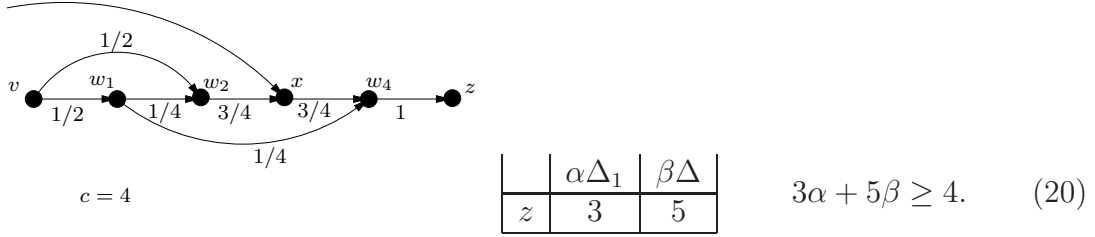
\begin{equation} \label{app:eq16}
\input redge2ciii1a.pstex_t
\hspace{0.5cm}
\begin{tabular}{|l|c|c|}
& $\alpha \Delta_1$ & $\beta \Delta$ \\
\hline
$z$ & 3 & 5 \\
\hline
\end{tabular}
\hspace{1cm}
3\alpha + 5\beta\ge 4 .
\end{equation}

Suppose next that the unique downward neighbor $y$ of $x$ is not $w_4$.
The local structure, table and inequality look like this 
($y$ is drawn above $w_4$ because this yields a sharper inequality):
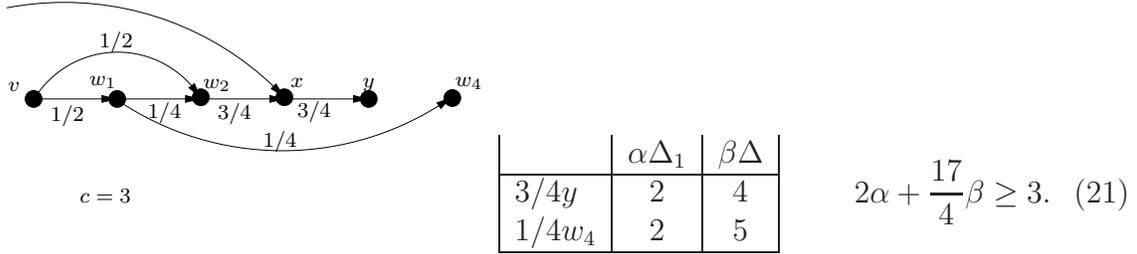
\begin{equation} \label{app:eq17}
\input redge2ciii1b.pstex_t
\hspace{0.5cm}
\begin{tabular}{|l|c|c|}
& $\alpha \Delta_1$ & $\beta \Delta$ \\
\hline
$3/4 y$ & 2 & 4 \\
$1/4 w_4$ & 2 & 5 \\
\hline
\end{tabular}
\hspace{1cm}
2\alpha + \frac{17}{4}\beta\ge 3 .
\end{equation}

\paragraph{Case 2.c.iii.2:} $x$ is a 2-vertex that precedes $w_4$.
This subcase splits into several subcases, where we assume,
respectively, that $\Delta(w_4)\ge 6$,
$\Delta(w_4) = 4$, and $\Delta(w_4) = 5$.

\paragraph{Case 2.c.iii.2(a).} 
Suppose first that $\Delta(w_4)\ge 6$. The configuration
looks like this:
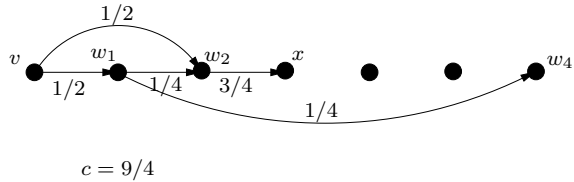
\begin{center}
\input redge2ciii2.pstex_t
\end{center}
The table and inequality are
\begin{equation} \label{app:eq18}
\begin{tabular}{|l|c|c|}
& $\alpha \Delta_1$ & $\beta \Delta$ \\
\hline
$3/4 x$ & 1 & 3 \\
$1/4 w_4$ & 1 & 6 \\
\hline
\end{tabular}
\hspace{2cm}
\alpha + \frac{15}{4}\beta\ge \frac{9}{4} .
\end{equation}

\paragraph{Case 2.c.iii.2(b).} 
Suppose next that $\Delta(w_4)=4$, and that one of the downward neighbors 
of $x$ is $w_4$.
Let $z$ denote the other downward neighbor. $w_4$ is a 1-vertex, and we
denote by $w$ its unique downward neighbor.
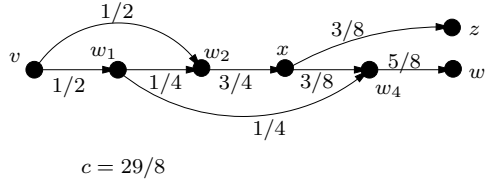
\begin{center}
\input redge2ciii2a.pstex_t
\end{center}
The 3-connectivity of the edge graph of $P$ implies, as above, that
$w\ne z$. 
Since we assume that $\Delta(w_4)=4$, $z$ 
also lies below $w_4$, and the table and inequality are
\begin{equation} \label{app:eq19}
\begin{tabular}{|l|c|c|}
& $\alpha \Delta_1$ & $\beta \Delta$ \\
\hline
$5/8 w$ & 2 & 5 \\
$3/8 z$ & 2 & 6 \\
\hline
\end{tabular}
\hspace{2cm}
2\alpha + \frac{43}{8}\beta\ge \frac{29}{8} .
\end{equation}

Suppose next that $\Delta(w_4)=4$ and $w_4$ is not a downward 
neighbor of $x$. Denote those two neighbors as $w$ and $z$, both of
which lie lower than $w_4$, by assumption, and are clearly
distinct. The configuration, table and inequality look like this:
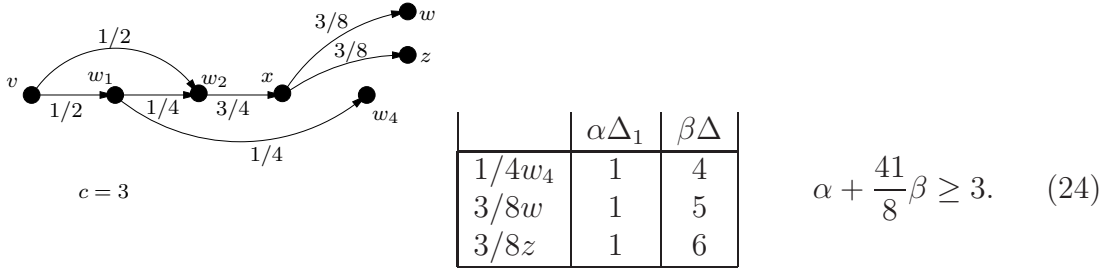
\begin{equation} \label{app:eq20}
\input redge2ciii2b.pstex_t
\hspace{0.5cm}
\begin{tabular}{|l|c|c|}
& $\alpha \Delta_1$ & $\beta \Delta$ \\
\hline
$1/4 w_4$ & 1 & 4 \\
$3/8 w$ & 1 & 5 \\
$3/8 z$ & 1 & 6 \\
\hline
\end{tabular}
\hspace{1cm}
\alpha + \frac{41}{8}\beta\ge 3 .
\end{equation}

\paragraph{Case 2.c.iii.2(c).} 
It remains to consider the case $\Delta(w_4)=5$. Let $z$ denote the 
unique vertex lying between $x$ and $w_4$. We may assume that $z$ is 
connected to $x$, for otherwise $z$ is not reachable from $v$, and 
we might as well reduce this case to the case $\Delta(w_4)=4$ just treated.

Consider first the subcase where the other downward neighbor of
$x$ is $w_4$ itself. Then $w_4$ is a 1-vertex, and we denote by
$w$ its unique downward neighbor. This subcase splits further into 
two subcases: First, assume that $z$ is a 1-vertex, and let $y$
denote its unique downward neighbor. Clearly, $y$ must lie below
$w_4$ (it may coincide with or precede $w$). The configuration
looks like this:
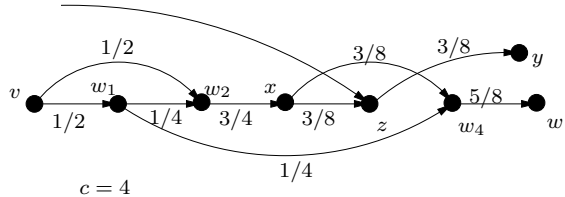
\begin{center}
\input redge2civ1a.pstex_t
\end{center}
The table and inequality are
\begin{equation} \label{app:eq21}
\begin{tabular}{|l|c|c|}
& $\alpha \Delta_1$ & $\beta \Delta$ \\
\hline
$3/8 y$ & 3 & 6 \\
$5/8 w$ & 3 & 6 \\
\hline
\end{tabular}
\hspace{2cm}
3\alpha + 6\beta\ge 4 .
\end{equation}

In the other subcase, $z$ is a 2-vertex; we denote its two
downward neighbors as $y$ and $t$. The vertices $w,y,t$ all lie
below $w_4$ and may appear there in any order. The configuration
looks like this:
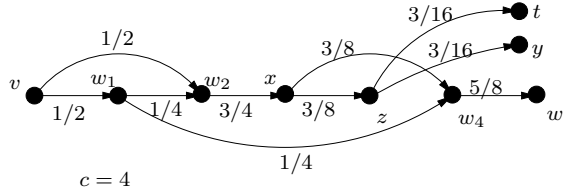
\begin{center}
\input redge2civ1b.pstex_t
\end{center}
The table and inequality are 
\begin{equation} \label{app:eq22}
\begin{tabular}{|l|c|c|}
& $\alpha \Delta_1$ & $\beta \Delta$ \\
\hline
$3/16 y$ & 2 & 6 \\
$3/16 t$ & 2 & 7 \\
$5/8 w$ & 2 & 6 \\
\hline
\end{tabular}
\hspace{2cm}
2\alpha + \frac{99}{16}\beta\ge 4 .
\end{equation}

Consider next the subcase where $w_4$ is not a downward neighbor of
$x$. Denote the other downward neighbor of $x$ as $y$, which lies
strictly below $w_4$.
This subcase splits into three subcases. First, assume that
$z$ is a 1-vertex, and denote its unique downward neighbor as $w$.
The configuration looks like this:
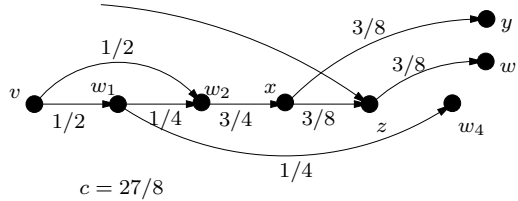
\begin{center}
\input redge2civ2a.pstex_t
\end{center}
The table and inequality are
\begin{equation} \label{app:eq23}
\begin{tabular}{|l|c|c|}
& $\alpha \Delta_1$ & $\beta \Delta$ \\
\hline
$1/4 w_4$ & 2 & 5 \\
$3/8 y$ & 2 & 6 \\
$3/8 w$ & 2 & 5 \\
\hline
\end{tabular}
\hspace{2cm}
2\alpha + \frac{43}{8}\beta\ge \frac{27}{8} .
\end{equation}

Second, assume that $z$ is a 2-vertex, so that none of its two
downward neighbors is $w_4$. Denote these neighbors as $w$ and
$t$. All three vertices $y,t,w$ lie strictly below $w_4$.
The configuration looks like this:
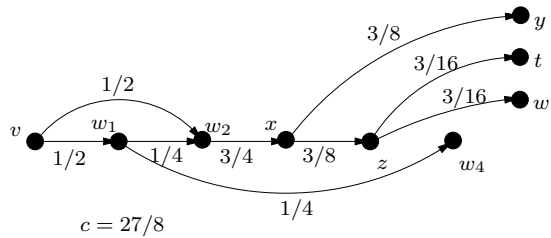
\begin{center}
\input redge2civ2b.pstex_t
\end{center}
The table and inequality are
\begin{equation} \label{app:eq24}
\begin{tabular}{|l|c|c|}
& $\alpha \Delta_1$ & $\beta \Delta$ \\
\hline
$1/4 w_4$ & 1 & 5 \\
$3/8 y$ & 1 & 6 \\
$3/16 w$ & 1 & 6 \\
$3/16 t$ & 1 & 7 \\
\hline
\end{tabular}
\hspace{2cm}
\alpha + \frac{95}{16}\beta\ge \frac{27}{8} .
\end{equation}

Finally, assume that $z$ is a 2-vertex, so that one of its two
downward neighbors is $w_4$. Denote the other neighbor as $w$. In
this case $w_4$ is a 1-vertex, and we denote its unique downward
neighbor as $t$. All three vertices $y,t,w$ lie strictly below $w_4$.
The configuration looks like this:
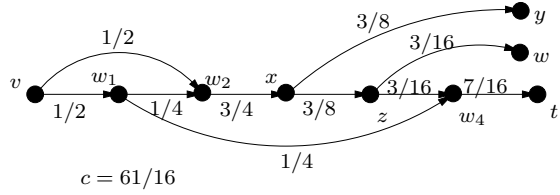
\begin{center}
\input redge2civ2c.pstex_t
\end{center}
The table and inequality are
\begin{equation} \label{app:eq25}
\begin{tabular}{|l|c|c|}
& $\alpha \Delta_1$ & $\beta \Delta$ \\
\hline
$3/8 y$ & 2 & 6 \\
$3/16 w$ & 2 & 7 \\
$7/16 t$ & 2 & 6 \\
\hline
\end{tabular}
\hspace{2cm}
2\alpha + \frac{99}{16}\beta\ge \frac{61}{16} ,
\end{equation}
which, by the way, is weaker than (\ref{app:eq22}).


\end{document}